\documentclass[12pt]{amsart}
\usepackage{color}
\usepackage{epsfig}
\usepackage{verbatim}

\begin{document}
\title {Splitters and decomposers for binary matroids}
\maketitle 
 
\begin {center}
S. R. Kingan 
\footnote{The author is partially supported by  PSC-CUNY grant number 66305-00 44.} \\     
Department of Mathematics \\
Brooklyn College, City University of New York\\
 Brooklyn, NY 11210\\
skingan@brooklyn.cuny.edu\\  
\end {center}
\bigskip 

\begin{abstract} Let $EX[M_1\dots, M_k]$ denote the class of binary matroids with no minors isomorphic to $M_1, \dots, M_k$.  In this paper we give a decomposition theorem for  $EX[S_{10}, S_{10}^*]$, where $S_{10}$ is a certain 10-element rank-4 matroid. As corollaries we obtain decomposition theorems for the classes obtained by excluding the Kuratowski graphs $EX[M(K_{3,3}), M^*(K_{3,3}), M(K_5), M^*(K_5)]$ and $EX[M(K_{3,3}), M^*(K_{3,3})]$. These decomposition theorems imply results on internally $4$-connected matroids by Zhou [\ref{Zhou2004}], Qin and Zhou [\ref{Qin2004}], and Mayhew, Royle and Whittle [\ref{Mayhewsubmitted}].
\end{abstract}

\bigskip 

\section {\bf Introduction}

Some matroids such as the complete graph on five vertices $K_5$, the complete bipartite graph with three vertices in each class $K_{3,3}$, and the Fano matroid $F_7$ play a more prominent role in structure theory than others.  Let $EX[M_1\dots, M_k]$ denote the class of binary matroids with no minors isomorphic to $M_1, \dots, M_k$.  In 1937 Wagner building on Kuratowski's work identified the minimal excluded minors for the class of planar graphs as $K_5$ and $K_{3,3}$ [\ref{Oxley2012}, 5.2.5]. In 1958 Tutte identified the minimal excluded minors for the  class of regular matroids as $F_7$ and its dual $F_7^*$ [\ref{Oxley2012}, 10.1.2].  
In 1980 Seymour gave a decomposition result for the class of regular matroids. In his decomposition result, we meet two matroids $R_{10}$ and $R_{12}$ that play a central role for regular matroids. The matroid $R_{10}$ is a 10-element rank-5 self-dual matroid. It is a splitter for 3-connected regular matroids. The matroid $R_{12}$ is a 12-element rank-6 self-dual matroid. It is a 3-decomposer for 3-connected regular matroids. Seymour proved that if $M$ is a 3-connected matroid in $EX[F_7, F_7^*]$, then either $R_{12}$ is a 3-decomposer for $M$ or $M$ is a graph, cograph, or $R_{10}$ [\ref{Oxley2012}, 13.1.2]. We talk of 3-connected matroids because matroids that are not 3-connected can be decomposed into their 3-connected components by direct sums or 2-sums [\ref{Oxley2012}, 8.3.1]. 

From a structural point of view, the regular matroids are taken care of and the focus shifts to matroids that have $F_7$ or $F_7^*$ as a minor. Observe that $F_7$ is the binary projective plane $PG(2,2)$ and has no further 3-connected single-element extensions. Matrix representations for $F_7$ and its dual $F_7^*$ are given below. By a slight abuse of notation, we will treat a vector matroid $M[A]$ as synonymous with the matrix $A$, when the context and representation is clear.  

\[ 
F_7=\left[ 
\begin{array}{cc|cccccc}
&&   0&1&1&1  \\
&I_3&1&0&1&1 \\
&&   1&1&0&1 
\end{array} 
\right] 
F_7*=\left[ 
\begin{array}{cc|ccccc}
&&   0&1&1  \\
&I_4&1&0&1 \\
&&   1&1&0 \\
&&   1&1&1
\end{array} 
\right] 
\]

\noindent The 3-connected single-element extensions of $F_7^*$ are $AG(3,2)$ and $S_8$, both of which are self-dual, so they are also the 3-connected coextensions of $F_7$. 

\[ 
AG(3,2)=\left[ 
\begin{array}{cc|cccccc}
&&   0&1&1&1   \\
&I_4&1&0&1&1  \\
&&   1&1&0&1  \\
&&   1&1&1&0 
\end{array} 
\right] 
S_8=\left[ 
\begin{array}{cc|ccccc}
&&   0&1&1&1   \\
&I_4&1&0&1&1  \\
&&   1&1&0&1  \\
&&   1&1&1&1 
\end{array} 
\right] 
\] 
\normalsize

$S_8$ has two 3-connected binary non-isomorphic single-element extensions $P_9$ and $Z_4$ and $AG(3,2)$ has only one single-element extension $Z_4$. $P_9$ has three single-element extensions, one of which is the internally 4-connected matroid $S_{10}$. In [\ref{Oxley1987}] Oxley flagged $P_9$ as important by obtaining a complete characterization of the 3-connected binary non-regular matroids with no minor isomorphic to $P_9$ or $P_9^*$. This made the new starting point for investigating binary matroids $P_9$ or $P_9^*$.  Observe that $P_9$ has a non-minimal exact 3-separation $(A, B)$, where $A=\{1,2,5,6\}$ is both a circuit and a cocircuit. Matrix representations for $P_9$ and $S_{10}$ are given below:

\[ 
P_9=\left[ 
\begin{array}{cc|ccccc}
&&   0&1&1&1&1  \\
&I_4&1&0&1&1&1 \\
&&   1&1&0&1&0 \\
&&   1&1&1&1&0
\end{array} 
\right]  
S_{10}=\left[ 
\begin{array}{c|cccccc}
&   0&1&1&1&1&1 \\
I_4&1&0&1&1&1&0 \\
&   1&1&0&1&0&0 \\
&   1&1&1&1&0&1
\end{array} 
\right] 
\]

In this paper we highlight $S_{10}$ and give a decomposition   theorem for $EX[S_{10}, S_{10}^*]$.   In the process we flag two of the single-element coextensions of $P_9$ known as $E_4$ and $E_5$ as significant. Matrix representations are shown below:

\[ 
E_4=\left[ 
\begin{array}{c|ccccc}
&    0&1&1&1&1 \\
&    1&0&1&1&1 \\
I_5& 1&1&0&1&0 \\
&    1&1&1&1&0 \\
&    0&1&0&0&1
\end{array} 
\right]  
E_5=\left[ 
\begin{array}{c|ccccc}
&    0&1&1&1&1 \\
&    1&0&1&1&1 \\
I_5& 1&1&0&1&0 \\
&    1&1&1&1&0 \\
&    1&0&1&0&0
\end{array} 
\right] 
\]

The matroids $E_5$ is internally 4-connected, whereas $E_4$ is not internally 4-connected. All these matroids play signficiant roles in the characterization of the almost-graphic matroids and the almost-regular matroids with at least one regular element [\ref{Kingan2002}]. The matroid $S_{10}$ is the first matroid in the internally 4-connected infinite family of almost-graphic matroids $S_{3n+1}$. The almost-regular matroids with at least one regular element are precisely the almost-regular matroids with no $E_5$-minor. Finally, the  matroid $T_{12}$ also makes an appearance in our characterization. It is a 12-element rank-6 self-dual 4-connected matroid [\ref{Kingan1996}]. A matrix representation is given below.

\[  
T_{12}=\left[ 
\begin{array}{c|cccccc}
&   1&1&0&0&0&1 \\
&   1&0&0&0&1&1 \\
I_6&0&0&0&1&1&1 \\
&   0&0&1&1&1&0 \\
&   0&1&1&1&0&0 \\
&   1&1&1&0&0&0 
\end{array} 
\right] 
\]

The next two theorems are the main results of this paper. In Theorem 1.1 we give a set of splitters and 3-decomposers that together characterize the excluded minor class $EX[S_{10}, S_{10}^*]$. Unlike the class of regular matroids that has just one 3-decomposer $R_{12}$ whose non-minimal exact 3-separation is induced in all 3-connected regular matroids containing it, $EX[S_{10}, S_{10}^*]$ has several 3-decomposers with different inducers. In fact, $E_4$ has two different non-minimal exact 3-separations and it is not possible to determine conclusively which of the two gets induced in the 3-connected matroids with an $E_4$-minor. We only show that at least one of the two 3-separations must be an inducer. Thus there is no way of determining how the 3-connected matroids in $EX[S_{10}, S_{10}^*]$ are pieced together across 3-sums, which is something Seymour was able to do for regular matroids. Nonetheless, when the focus shifts from a particular non-minimal exact 3-separation to a set of 3-decomposers with different non-minimal exact 3-separations, we get a decomposition theorem with a relatively short proof. Several recent theorems on internally 4-connected matroids follow immediately.

\bigskip
\noindent {\bf Theorem 1.1.} {\it Let $M$ be a $3$-connected non-regular matroid in $EX[S_{10}, S_{10}^*]$. Then  $S_8$, $P_9$, $P_9^*$, or $E_4$ (with either of two inducers) is a $3$-decomposer for $M$ or $M$ is isomorphic to $F_7$, $F_7^*$, $E_5$, $T_{12}\backslash e$, $T_{12}/e$,or  $T_{12}$.}
\bigskip

Once we establish that a matroid in the class has a non-minimal exact 3-separation induced by one of the listed 3-decomposers, we are not concerned with it any more. We know it can be decomposed into smaller parts. It is the other ones that won't admit a  decomposition that are of interest. The rank 3 extremal matroid (monarch) is $F_7$ and the rank 6 extremal matroid is $T_{12}$. The rank 4 and 5 extremal matroids are not important because they admit non-minimal exact 3-separation. 

Next, observe that $S_{10}$ is the only single-element extension of $M^*(K_{3,3})$ and $E_5$ is a non-regular single-element coextensions of $M^*(K_{3,3})$.  A decomposition result for $EX[M(K_{3,3}), M^*(K_{3, 3})]$ follows directly from Theorem 1.1 because by excluding $S_{10}$ or $S_{10}^*$ we are effectively excluding $M(K_{3,3})$ or $M^*(K_{3, 3})$.  No additional work has to be done to get the next major new decomposition theorem; just remove $E_5$ from Theorem 1.1.

\bigskip
\noindent {\bf Theorem 1.2.} {\it Let $M$ be a $3$-connected binary non-regular matroid in $EX[M(K_{3,3}), M^*(K_{3, 3})]$. Then $S_8$, $P_9$, $P_9^*$, or $E_4$ (with either of two inducers) is a $3$-decomposer for $M$ or $M$ is isomorphic to  $F_7$, $F_7^*$,  $T_{12}\backslash e$, $T_{12}/e$,or  $T_{12}$.}
\bigskip

When the extremal matroids are known all the internally 4-connected matroids are easily found because they are the minors of the extremal matroids. The main proofs of three recent results follow immediately as corollaries. The first corollary is the main theorem in [\ref{Zhou2004}, Theorem 1.1].  The reader will want to know that Zhou had different names for the matroids. In his paper $E_4=P_{10}$, $E_5=N_{10}$, and $S_{10}= \tilde{K_5}$. Our notation of $E_4$ and $E_5$ is from the almost-graphic paper and a dozen other papers. So we want to stick with it for consistency. Additionally, Zhou also talked about a single-element coextension of $S_{10}$ called $Q_{10}$ (which we call $E_7$), but we show in the proof of the main theorem that it is not significant because it is decomposed by $P_9$. 

\bigskip
\noindent {\bf Corollary  1.3.} {\it (Zhou 2004) Let $M$ be an internally $4$-connected binary non-regular matroid in $EX[S_{10}, S_{10}^*]$. Then $M$ is isomorphic to $F_7$, $F_7^*$, $E_5$, $T_{12}$, $T_{12}\backslash e$, or $T_{12}/e$.} $\qed$
\bigskip

We already know from [\ref{Kingan1997}] that $EX[M(K_{3,3}), M^*(K_{3,3}), M(K_5), M^*(K_5)]$ is mostly the same as $EX[M(K_{3,3}), M^*(K_{3,3})]$. The main theorem in [\ref{Kingan1997}, Theorem 2.1] states that, if $M$ is a $3$-connected binary matroid with an $M(K_5)$-minor, then either $M$ has an $M(K_{3,3})$- or $M^*(K_{3,3})$-minor or $M$ is isomorphic to $M(K_5)$, $T_{12}/e$, or $T_{12}$. Moreover,  $T_{12}$ is a splitter for $EX[M(K_{3,3}), M^*(K_{3,3})]$. It follows that matroids in  $EX[M(K_{3,3}), M^*(K_{3,3})]$, but not in $EX[M(K_{3,3}), M^*(K_{3,3}), M(K_5), M^*(K_5)]$ are precisely $M(K_5)$, $M^*(K_5)$, $T_{12}\backslash e$, $T_{12}/e$, or $T_{12}$.

 Mayhew, Royle, and Whittle's identification of the internally 4-connected matroids in \break $EX[M(K_{3,3}), M^*(K_{3, 3})]$ in [\ref{Mayhewsubmitted}] follows directly from Theorem 1.3.  So does the main result of Qin and Zhou's paper where they exclude the Kuratowski graphs and their duals from binary matroids [\ref{Qin2004}, Theorem 1.3].  Note that if a matroid $M$ in $EX[M(K_{3,3}), M^*(K_{3,3}), M(K_5), M^*(K_5)]$  is regular, then $M$ is isomorphic to a planar graph. If a matroid in $EX[M(K_{3,3}), M^*(K_{3,3})]$  is regular, then $M$ is isomorphic to a planar graph or $K_5$ or $K_5^*$.

\bigskip
\noindent {\bf Corollary  1.4.} {\it (Mayhew, Royle, Whittle) Let $M$ be an internally $4$-connected binary non-regular matroid in   $EX[K_{3,3}, K_{3, 3}^*]$. Then $M$ is isomorphic to $F_7$, $F_7^*$, $T_{12}$, $T_{12}\backslash e$, or $T_{12}/e$.} $\qed$
\bigskip

\noindent {\bf Corollary  1.5.} {\it (Qin and Zhou 2004) Let $M$ be an internally $4$-connected binary non-regular  matroid in   $EX[M(K_{3,3}), M^*(K_{3,3}), M(K_5), M^*(K_5)]$. Then $M$ is isomorphic to $F_7$ or  $F_7^*$.} $\qed$
\bigskip

In Section 2 we explain the terminology and techniques used in the paper. In Section 3 we give the proof of Theorem 1.1. Note that Theorem 1.2 and the corollaries are immediate consequences of Theorem 1.1.
 

\section {The techniques used in the main theorem}

The matroid terminology follows Oxley [\ref{Oxley2012}]. If $M$ and $N$ are matroids on the sets $E$ and $E \cup e$ where $e   \not\in E$, then $M$ is a  {\it single-element extension} of $N$ if $M \backslash e = N$,  and $M$ is a {\it single-element coextension} of $N$ if $M^*$ is a single-element extension of $N^*$. Let $\mathcal M$ be a class of matroids closed under minors and isomorphism.  A {\it splitter} $N$ for $\mathcal M$ is a 3-connected matroid in $\mathcal M$ such that no $3$-connected matroid in $\mathcal M$ has $N$ as a proper minor. Checking if a matroid is a splitter is a potentially infinite task. But the Splitter Theorem [\ref{Oxley2012}, Theorem 12.1.2] establishes that if every $3$-connected single-element extension and coextension of $N$ does not belong to $\mathcal M$, then $N$ is a splitter for $\mathcal M$.

Let $M$ be a matroid and $X$ be a subset of the ground set $E$. The {\it connectivity function}  $\lambda$ is defined as $\lambda (X) = r(X) + r(E-X) - r(M)$.  Observe that $\lambda (X) = \lambda (E-X)$.  For $k\ge 1$, a partition $(A, B)$ of $E$ is called a $k$-separation if $\lambda (A) \le k-1$ for $|A|, |B| \ge k$.  If $\lambda (A)=k-1$, we call $(A, B)$ an {\it exact k-separation}.  If $\lambda (A)=k-1$ and $|A|=k$ or $|B|=k$, we call $(A, B)$ a {\it minimal exact k-separation}. For $n\ge 2$, we say $M$ is {\it n-connected} if $M$ has no $k$-separation for $k\le n-1$.  A  $k$-connected matroid is {\it internally $(k+1)$-connected} if it has no non-minimal exact $k$-separations. 

Let $M$ be a matroid in  $\mathcal M$ with an $N$-minor and let $N$ have an exact  $k$-separation $(A, B)$. If there exists a $k$-separation $(X, Y)$ of $M$ such that $A\subseteq X$ and  $B\subseteq Y$, then we say the $k$-separation $(A, B)$ of $N$ is {\it induced} in $M$. Suppose $M$ is a $k$-connected matroid with a $k$-connected minor $N$ such that $N$ has a non-minimal exact $k$-separation $(A, B)$. We call $N$ a {\it $k$-decomposer} for $M$ having $(A, B)$ as an inducer, if $M$ has a non-minimal exact $k$-separation $(X, Y)$ such that $A\subseteq X$ and $B\subseteq Y$. 
If $(A, B)$ is not induced in $M$, then we say $M$ {\it bridges} the $k$-separation $(A, B)$ in $N$. Define $k_M(A, B)=min\{\lambda_M(X) | A\subseteq X \subseteq E(M)-B\}$. Thus $M$ bridges $(A, B)$ if and only if $k_M(A, B)\ge k$. 

The next theorem appears in [\ref{Kingansubmitted}, Theorem 1.2]. It gives a sufficient conditions to determine when an exact $k$-separation $(A, B)$ in $N$ is induced in $M$. Note that no isomorphism is involved in the calculations of  extensions and coextensions.
\bigskip

\noindent{\bf Theorem 2.1.} {\it  Let $N$ be a simple and cosimple matroid in $\mathcal{M}$ with an exact $k$-separation $(A,B)$, such that $A$ is the union of circuits and the union of cocircuits.  Suppose $M\in \mathcal {M}$.

\begin{enumerate}
\item[(i)] If $M$ is a simple single-element extension of $N$ such that $M\backslash e=N$ (no isomorphism is involved), then $\lambda_M(A)=k-1$ or $\lambda_M(A\cup e)=k-1$.

\item[(ii)] If $M$ is a cosimple single-element coextension of $N$ such that $M/ f=N$, then $\lambda_M(A)=k-1$ or $\lambda_M(A\cup f)=k-1$.

\item[(iii)] If $M$ is a cosimple single-element coextension of a Type (i) matroid or a simple single-element extension of a Type (ii) matroid, then $M$ satisfies one of the following conditions:

\begin{enumerate}
\item[(a)] $\lambda _{M/f}(A)=k-1$ and $\lambda _{M\backslash e}(A)=k-1$;
\item[(b)] If $\lambda _{M/f}(A)=k-1$ and $\lambda _{M\backslash e}(A\cup f)=k-1$, then $\lambda _M(A\cup f)=k-1$ or $\{e, f, g\}$ is a triad or triangle with $g\in A$;
\item[(c)]If $\lambda _{M/f}(A\cup e)=k-1$ and $\lambda _{M\backslash e}(A)=k-1$, then $\lambda _M(A\cup e)=k-1$ or $\{e, f, g\}$ is a triad or triangle with $g\in A$;  or 
\item[(d)] If $\lambda _{M/f}(A\cup e)=k-1$ and $\lambda _{M\backslash e}(A\cup f)=k-1$, then $\{e, f, g\}$ is a triangle or triad in $M$ with $g\in A$.
\end{enumerate}
\end{enumerate}
\noindent Then the $k$-separation $(A, B)$ of $N$ is induced in $M$ for every $M\in \mathcal M$ with $N$ as a minor.}

\bigskip

Practically speaking we are interested in 3-decomposers for classes of $GF(q)$-representable matroids closed under minors and isomorphism because we have a ``Splitter Theorem" only for 3-connected matroids. A simple matroid is 3-connected if $\lambda (A)\ge 2$ for all partitions $(A, B)$ with $|A|\ge 3$ and $|B|\ge 3$. If $N$ is 3-connected, a simple single-element extension and a cosimple single-element coextension are also 3-connected. A 3-connected matroid is {\it internally $4$-connected} if $\lambda (A)\ge 3$ for all partitions $(A, B)$ with $|A|\ge 4$ and $|B|\ge 4$. In this case $\lambda (A)= 2$ is allowed only when either $|A|$ or $|B|$ has size at most 3. 
Suppose $M$ is a 3-connected matroid having a 3-connected minor $N$ and $N$ has a non-minimal exact 3-separation $(A, B)$. If $N$ is a {\it 3-decomposer} for $M$, then $M$ has a $3$-separation $(X, Y)$ such that $A\subseteq X$ and $B\subseteq Y$. In this case $\lambda (X)=2$ and $|X|\ge 4$, $|Y|\ge 4$. Hence, if $N$ is a 3-decomposer for $M$, then $M$ is not internally 4-connected. The converse is not true.

Next, we discuss the situation when $N$ has two distinct non-minimal exact 3-separations $(A_1, B_1)$ and $(A_2, B_2)$ where both $A_1$ and $A_2$ are unions of circuits and unions of cocircuits and both $(A_i, B_i)$ satisfy Theorem 2.1(i, ii). If both satisfy Theorem 2.1(i, ii) in such a way that $\lambda(A_i)=2$, then a one-element check suffices and we can conclude that both $(A_1, B_1)$ and $(A_2, B_2)$ are induced in all matroids in $\mathcal M$ with an $N$-minor. 
Problems arise when $\lambda (A_i\cup e)=2$. Now, suppose further that $(A_i, B_i)$ do not satisfy Condition (iii) for some of the coextension rows that can be added to Type (i) matroids. Then we cannot conclude that 
$(A_i, B_i)$ are induced in all matroids in $\mathcal M$ with an $N$-minor. Note that, if $M$ is not a bridging matroid, Condition (iii) fails  because $\lambda (A\cup e\})\neq 2$ or $\lambda (A\cup f\})\neq 2$, as the case may be. It is still true that $\lambda (A\cup \{e, f\})=2$ (or else it would be a bridging matroid). 
Let $B_i$ be the set of ``bad" coextension rows that cause Condition (iii) to fail for $(A_i, B_i)$, where $i=\{1, 2\}$. As long as $B_1$ and $B_2$ are disjoint, we can conclude that either $(A_1, B_1)$ or $(A_2, B_2)$ is induced in $M$ without being able to specify which one exactly. 

The next corollary follows from Theorem 2.1 and the above discussion. Further we hypothesize that $N$ is self-dual so the argument for the simple single-element extension for Type (ii) matroids follows by duality.
\bigskip

\noindent{\bf Corollary 2.2.} {\it  Let $N$ be a self-dual simple and cosimple matroid in $\mathcal{M}$ with two non-minimal exact $3$-separations $(A_1,B_1)$ and $(A_2, B_2)$, such that each of $A_1$ and $A_2$ are the union of circuits and the union of cocircuits.  Suppose $M\in \mathcal {M}$ is a Type (i) or (ii) matroid and condition (i) and (ii) hold for both $(A_1,B_1)$ and $(A_2, B_2)$ in such a way that if $\lambda(A_i\cup x)=2$, then $\lambda(A_j)=2$ for $i, j\in \{1, 2\}$ and $x\in \{e, f\}$. If $M$ is a cosimple single-element coextension of a Type (i) matroid such that condition (iii) is satisfied by either $(A_1, B_1)$ or $(A_2, B_2)$, then either $(A_1, B_1)$ or $(A_2, B_2)$ is induced in $M$ for every $M\in \mathcal M$ with $N$ as a minor. } $\qed$

\bigskip

We end this section by describing our method for calculating extensions and coextensions. Let $N$ be a $GF(q)$-representable $n$-element rank-$r$  matroid  represented by the matrix $A=[I_r|D]$ over $GF(q)$. The columns of A may be viewed as a subset of the columns of the matrix that represents the projective geometry $PG(r - 1, q)$. Let $M$ be a simple single-element extension of $N$ over $GF(q)$. Then $N=M\backslash e$ and $M$ may be represented by $[I_r|D']$, where $D'$ is the same as $D$, but with one additional column corresponding to the element $e$. The new column is distinct from the existing columns and has at least two non-zero elements. If the existing columns are labeled $\{1, \dots , r, \dots , n\}$, then the new column is labeled $(n+1)$.

Suppose $M$ is a cosimple single-element coextension of $N$ over $GF(q)$. Then $N=M/f$ and $M$ may be represented by the matrix $[I_{r+1}| D'']$, where $D''$ is the same as $D$, but with one additional row. The new row is distinct from the existing rows and has at least two non-zero elements. The columns of $[I_{r+1}| D'']$ are labeled $\{1, \dots , r+1, r+2, \dots , n, n+1\}$. The coextension element $f$ corresponds to column $r+1$. The coextension row is selected from $PG(n-r, q)$, which means there could be a much larger selection of row vectors for the coextension. We can visualize the new element $f$  as appearing in the new dimension and lifting several points into the higher dimension. Observe that $f$ forms a cocircuit with the elements corresponding to the non-zero elements in the new row.  Note that in $[I_{r+1}|D'']$ the labels of columns beyond $r$ are increased by 1 to accomodate the new column $r+1$. 

We refer to the simple single-element extensions of $N$ as Type (i) matroids and the cosimple single-element coextensions of $N$ as Type (ii) matroids. The structure of type (i) and Type (ii) matroids are shown in Figure 1.

\begin{figure}[h]
\centering
\epsfxsize 5in \epsfbox{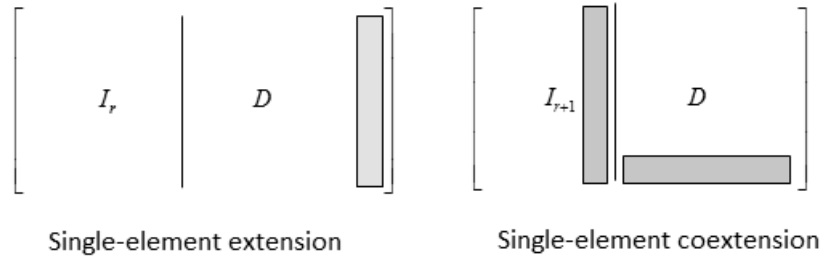}
\caption{Structure of Type (i) and Type (ii) matroids }
\end{figure}

Once the simple single-element extensions (Type (i) matroids) and cosimple single-element coextensions (Type (ii) matroids) are determined, the number of permissable rows and columns give a bound on the choices for the cosimple single-element extensions of the Type (i) matroids and the simple single-element extensions of the Type (ii) matroids, respectively. 

The structure of the cosimple single-element coextensions of a Type (i) matroid and the simple single-element extensions of a Type (ii) matroid are shown in Figure 2.
 
\begin{figure}[h]
\centering
\epsfxsize 5in \epsfbox{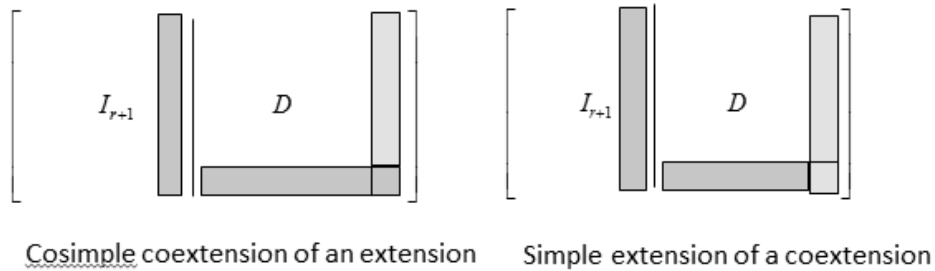}
\caption{Structure of $M$, where $|E(M)-E(N)|=2$ }
\end{figure}
 
\noindent When computing the cosimple single-element coextension of a Type (i) matroid, there are three types of rows that may be inserted into the last row.

\begin {enumerate}
\item[(i)] rows that can be added to $N$ to obtain a coextension with a 0 or 1 as the last entry (or as many as the entries in $GF(q)$ for higher order fields); 
\item[(ii)] the identity rows with a 1 in the last position; and
 \item[(iii)] rows ``in-series" to the right-hand side of the matrix with the last entry reversed. 
\end{enumerate}

\noindent When computing the simple single-element extension of a Type (ii) matroid, there are three types of rows that may be inserted into the last column.

\begin {enumerate}
\item[(i)] columns that can be added to $N$ to obtain an extension with a 0 or 1 as the last entry (or as many as the entries in $GF(q)$ for higher order fields); 
\item[(ii)] the identity columns with a 1 in the last position; and
 \item[(iii)] columns ``in-parallel" to the right-hand side of matrix with the last entry reversed. 
\end{enumerate}

\bigskip

\section {Proof of the main theorem}
 
\noindent {\bf Proof of Theorem 1.1.} Let $M$ be a $3$-connected non-regular matroid in $EX[S_{10}, S_{10}^*]$. Then $M$ has an $F_7^*$-minor. To get the 3-connected single-element extensions of $F_7^*$, compare the columns in $F_7^*$ to the columns in  $PG(3, 2)$ shown below.

\[ 
PG(3,2)=\left[ 
\begin{array}{cc|ccccccccccc}
&&   0&0&0&0&1&1&1&1&1&1&1  \\
&I_4&0&1&1&1&0&0&0&1&1&1&1  \\
&&   1&0&1&1&0&1&1&0&0&1&1  \\
&&   1&1&0&1&1&0&1&0&1&0&1
\end{array} 
\right] 
\]

\noindent Add missing columns one by one and group them into isomorphism classes.  This gives us the two non-isomorphic single-element extensions $AG(3,2)$ and $S_8$ shown below. Only one column from $PG(3,2)$ (namely $[1 1 1 0]$)  may be added to $F_7^*$ to obtain $AG(3,2)$, whereas adding any of the five columns $[0 0 1 1]$, $[0 1 0 1]$, $[0 1 1 0]$, $[1 0 0 1]$ $[1 1 0 0 ]$ or $[1 1 1 1]$ gives a representation for $S_8$ [\ref{Oxley2012}, 12.2.4].

\medskip
\noindent {\bf Claim 1.} {\it $S_8$ is a $3$-decomposer for $EX[P_9, P_9^*]$}. 
\medskip

\noindent {\bf Proof.}  $S_8$ has a non-minimal exact 3-separation $(A, B)$ where $A=\{1, 2, 5, 6\}$ and $\lambda (A)=2$.  Observe that, $S_8$ has only two non-isomorphic single-element extensions $P_9$ and $Z_4$. Moreover, only one column $[1 1 1 0]$ may be added to obtain $Z_4$ and the remaining columns from $PG(3, 2)$ give $P_9$. It is easy to check that $\lambda (\{1, 2, 5, 6\})=2$ for this matrix representing $Z_4$. Since $S_8$ is self-dual, $Z_4^*$ is the only single-element coextension in $EX[P_9, P_9^*]$ and it is obtained by adding only one row $[1 1 1 0]$. We can check that $\lambda (\{1, 2, 6, 7\})=2$ in this coextension. By Theorem 2.1, $S_8$ is a 3-decomposer for $EX[P_9, P_9^*]$. 

\medskip
\noindent {\bf Claim 2.} {\it $P_9$ or $P_9^*$ are $3$-decomposers for $EX[S_{10}, S_{10}^*, E_4, E_5]$. Moreover, $E_5$ is internally $4$-connected.}
\medskip

\noindent {\bf Proof.} Observe that $P_9$ has three non-minimal exact 3-separations. They are $(A, B)$, $(A_1, B_1)$, and $(A_2, B_2)$, where $A=\{1, 2, 5, 6\}$, $A_1=\{3, 4, 7, 8\}$ and $A_2=\{3, 4, 7, 9\}$. While $A$ is a circuit and a cocircuit  $B$, $B_1$, and $B_2$ are not  unions of cocircuits and $A_1$ and $A_2$ are not unions of circuits. Thus the only candidate for a non-minimal exact 3-separation is $(A, B)$. Consider the single-element extensions of $P_9$:
\begin {itemize}
\item adding column $[1 1 1 0]$  gives $D_1$; 
\item adding any one of columns   $[1 0 0 1]$, $[0 1 0 1]$, $[0 1 1 0]$, or $[1 0 1 0]$  gives $S_{10}$; and 
\item adding column $[0 0 1 1]$    gives  $D_3$. 
\end{itemize}
For the first and third extension, $\lambda (\{1, 2, 5, 6\})=2$. $P_9$ has 8 non-isomorphic cosimple single-element coextensions. They are obtained by adding the following columns: 
\begin {itemize}
\item  $[1 1 0 0 0]$ or $[1 1 1 1 1]$  gives $E_1$;
\item  $[1 1 0 1 1]$ or  $[1 1 1 0 0]$    gives $E_2$; 
\item  $[1 1 0 0 1]$ or $[1 1 1 0 1]$    gives  $E_3$; 
\item  $[0 1 0 0 1]$, $[0 1 0 1 0]$, $[0 1 1 0 1]$, $[0 1 1 1 0]$, $[1 0 0 0 1]$,     $[1 0 0 1 0]$, $[1 0 1 0 1]$, or $[1 0 1 1 0]$    gives  $E_4$;
\item  $[0 1 0 1 1]$, $[0 1 1 0 0]$, $[1 0 0 1 1]$, or $[1 0 1 0 0]$   gives  $E_5$; 
\item  $[0 0 1 0 1]$ or $[0 0 1 1 0]$    gives  $E_6$;  
\item  $[0 0 1 1 1]$                  gives  $E_6^*$; and
\item $[0 0 0 1 1]$      gives  $E_7$. 
\end{itemize}
By checking each row, regardless of isomorphism, we see that $\lambda(\{1, 2, 6, 7\})=2$ for all the rows that give $E_1$, $E_2$, $E_3$, $E_6$, $E_6^*$, and $E_7$. The result follows from Theorem 2.1. Further note that  $E_4$ and $E_5$ are self-dual and $E_5$ is internally 4-connected.   
\medskip

\medskip
\noindent {\bf Claim 3.} {\it $E_5$ is a splitter for $EX[S_{10}, S_{10}^*]$. }
\medskip

\noindent {\bf Proof.} It is straightforward to check that every single-element extension of $E_5$ has a minor isomorphic to $S_{10}$.  The single-element extensions are obtained by adding the following columns:
\begin{itemize}
\item $[0 0 0 1 1]$, $[0 0 1 0 1]$, $[1 0 0 1 0]$, $[1 0 1 0 0]$ gives $(E_5, ext1)$;  
\item $[0 0 1 1 0]$ or $[1 0 0 0 1]$ gives $(E_5, ext2)$;
\item $[0 0 1 1 1]$, $[1 0 0 1 1]$, $[1 0 1 0 1]$ or $[1 0 1 1 0]$ gives $(E_5, ext3)$; 
\item $[0 1 0 0 1 ]$, $[0 1 1 0 0]$, $[0 1 1 1 1]$ or  $[1 1 1 0 1]$  gives $(E_5, ext4)$;
\item $[0 1 0 1 0]$, $[1 1 0 0 0]$, $[1 1 0 1 1]$ or $[1 1 1 1 0]$ gives $(E_5, ext5)$;
\item $[0 1 0 1 1]$ or $[1 1 1 0 0]$     gives  $(E_5, ext6)$; and
\item $[0 1 1 0 1]$  gives  $(E_5, ext7)$.
\end{itemize}
Note that since we are checking for minors it is sufficient to check one column per isomorphism class and $E_5$ is self-dual, every single-element coextension has a minor isomorphic to $S_{10}^*$. It follows that $E_5$ is a splitter for $EX[S_{10}, S_{10}^*]$. 
\medskip

Returning to the proof of Theorem 1.1 consider the simple single-element extensions and cosimple single-element coextensions of $E_4$. Observe that
\begin{itemize}
  \item adding column $[0 0 1 1 0]$ $[1 0 1 1 0]$ gives $A$; 
  \item adding column $[0 1 1 1 1]$ $[1 1 1 0 0]$ gives $B$; 
  \item adding column $[1 1 0 0 0]$ gives $C$; and 
  \item adding column $[1 1 0 1 1]$ gives $T_{12}/ e$.
\end{itemize}
Adding any other column to $E_4$ gives a matroid with an $S_{10}$-minor. Similarly, 
\begin{itemize}
  \item adding row $[0 0 1 1 0]$ or $[1 0 0 0 1]$ gives $A^*$; 
  \item adding row $[1 1 0 0 1]$ $[1 1 1 0 0]$ gives $B^*$; 
  \item adding row $[1 1 0 0 0]$   gives $C^*$; and  
  \item adding row $[0 1 0 1 0]$ gives $T_{12}\backslash e$
\end{itemize}
Adding any other row to $E_4$ gives a matroid with an $S_{10}$-minor:
From $T_{12}/e$ we get $T_{12}$ which is a splitter for $EX[S_{10}, S_{10}^*]$. See [\ref{Kingan1997}] for the details. Thus if $M$ has a $T_{12}/e$ as a minor, then $M\cong T_{12}/e, T_{12}\backslash e, T_{12}$.  

Next, observe that $E_4$ has two non-minimal exact 3-separations $(A_1, B_1)$ and $(A_2, B_2)$, where $A_1=\{1, 2, 5, 6, 7, 10\}$ and $A_2=\{1, 2, 3, 4, 8, 9\}$. Moreover, $A_1$ and $A_2$ are both union of circuits and union of cocircuits. Note that $B_1$ and $B_2$ do not meet this condition.

\smallskip
\centerline { $A_1=\{6, 7, 10\}\cup \{1, 2, 5, 10\}$ and $A_1=\{5, 7, 10\}\cup \{1, 2, 6, 10\}$}

\centerline { $A_2=\{3, 8, 9\}\cup \{1, 2, 4, 8\}$ and $A_2=\{3, 4, 8\}\cup \{1, 2, 3, 9\}$}

\smallskip

\noindent Using Theorem 2.1 and Corollary 2.2 we will prove that every matroid in $EX[S_{10}, S_{10}^*]$ with an $E_4$-minor has a non-minimal exact 3-separation induced by either $(A_1, B_1)$ or $(A_2, B_2)$. As we see in the proof, we cannot determine which of the two 3-separations is induced in a particular matroid, just that at least one of them is induced.

\medskip
\noindent {\bf Claim 4.} {\it Suppose $M\in EX[S_{10}, S_{10}^*]$ has an $E_4$-minor. Then $M$ has a non-minimal exact $3$-separation induced by either $(A_1, B_1)$ or $(A_2, B_2)$.}
\medskip

\noindent {\bf Proof.} Table 1a and 1b show that Theorem 2.1(i, ii) are satisfied.
 
\small
 \begin{center}
\begin{tabular}{|c|c|c|c|}\hline

\bf Name  & \bf Extension Columns   &  $ \bf A_1=\{1, 2, 5, 6, 7, 10\}$ & $ \bf A_2=\{1, 2, 3, 4, 8, 9\}$\\  \hline \hline

$A_{22}$ & $\alpha=[0 0 1 1 0]$  &  $\lambda\{1, 2, 5, 6, 7, 10\}=2$ & $\lambda\{1, 2, 3, 4, 8, 9, {\bf 11}\}=2$    \\  \hline

& $\beta = [1 0 1 1 0]$  &  $\lambda\{1, 2, 5, 6, 7, 10, {\bf 11}\}=2$ &  $\lambda\{1, 2, 3, 4, 8, 9\}=2$ \\  \hline \hline

$A_{11}$ & $\gamma = [0 1 1 1 1]$   &   $\lambda\{1, 2, 5, 6, 7, 10, {\bf 11}\}=2$  &  $\lambda\{1, 2, 3, 4, 8, 9\}=2$ \\   \hline
   
& $\delta = [1 1 1 0 0]$   &   $\lambda\{1, 2, 5, 6, 7, 10\}=2$ &  $\lambda\{1, 2, 3, 4, 8, 9, {\bf 11}\}=2$\\   \hline   \hline

$A_5$ & $\epsilon = [1 1 0 0 0]$   &   $\lambda\{1, 2, 5, 6, 7, 10\}=2$ & $\lambda\{1, 2, 3, 4, 8, 9\}=2$\\   \hline \hline

\end{tabular}
\end{center}
\smallskip
 \normalsize
 \begin{center} Table 1a: Simple single-element extensions of $E_4$ in $EX[S_{10}, S_{10}^*]$ \end{center} 

\smallskip
\small
 \begin{center}
\begin{tabular}{|c|c|c|c|}\hline

\bf Name  & \bf Coextension Rows   & $ \bf A_1=\{1, 2, 5, 7, 8, 11\}$ & $ \bf A_2=\{1, 2, 3, 4, 9, 10\}$\\  \hline \hline

$A_{22}^*$ & $a=[0 0 1 1 0]$  &  $\lambda\{1, 2, 5, 7, 8, 11\}=2$ &  $\lambda\{1, 2, 3, 4, {\bf 6}, 9, 10\}=2$  \\  \hline
        & $b=[1 0 0 0 1]$  &  $\lambda\{1, 2, 5, {\bf 6}, 7, 8, 11\}$ &  $\lambda\{1, 2, 3, 4, 9, 10\}=2$  \\  \hline \hline

$A_{11}^*$ & $c=[1 1 0 0 1]$   &   $\lambda\{1, 2, 5, {\bf 6}, 7, 8, 11\}=2$ &  $\lambda\{1, 2, 3, 4, 9, 10\}=2$\\   \hline
   
        & $d=[1 1 1 0 0]$   &   $\lambda\{1, 2, 5, 7, 8, 11\}=2$ &  $\lambda\{1, 2, 3, 4, {\bf 6}, 9, 10\}=2$\\   \hline   \hline

$A_5^*$ & $e=[1 1 0 0 0]$   &   $\lambda\{1, 2, 5, 7, 8, 11\}=2$ &  $\lambda\{1, 2, 3, 4, 9, 10\}=2$\\   \hline \hline

\end{tabular}
\end{center}
\smallskip
\normalsize
\begin{center} Table 1b: Cosimple single-element coextensions of $E_4$ in $EX[S_{10}, S_{10}^*]$\end{center} 
\smallskip

Next, we compute the cosimple single-element coextensions of the Type (i) matroids $A$, $B$, and $C$. As explained in Section 2, there are three types of rows that can be added: rows $a$ to $e$ with a 0 or a one in the last entry; the identity rows with a 1 in the last entry; and rows in-series with rows in the matrix with the last entry switched. We handle $(A_1, B_1)$ and $(A_2, B_2)$ separately. However, notice that $A$ with column $[1 0 1 1 0]$ and $B$ with column $[0 1 1 1 1]$ have coextensions that exhibit similar behavior.

\smallskip
\small
 \begin{center}
\begin{tabular}{|c|p{10em}|p{6em}|c|p{13em}|}
\hline
\bf   &\bf{Coextension Rows} &\bf{$S_{10}$, $S^*_{10}$ } & \bf {$A_1=\{1, 2, 5, 7, 8, 11\}$} \\  \hline \hline

$A$ with column $[1 0 1 1 0]$  & $a=[0 0 1 1 0 0]$  & No &   $\lambda\{1, 2, 5, 7, 8, 11, {\bf 12}\}=2$\\ \hline
and & $a'=[0 0 1 1 0 1]$ 	& Yes 	& -           \\  \hline
$B$ with column $[0 1 1 1 1]$  & $b=[1 0 0 0 1 0]$ 	& No 		&  Bad row		\\  \hline
& $b'=[1 0 0 0 1 1]$ 	& No 		&  Bad row		 \\  \hline
& $c=[1 1 0 0 1 0]$ 	& No 		& Bad row 		\\  \hline
& $c'=[1 1 0 0 1 1]$ & No  	& Bad row   		\\  \hline
& $d=[1 1 1 0 0 0]$  & Yes 	& -          \\  \hline  
& $d'=[1 1 1 0 0 1]$ & No 		& $\lambda\{1, 2, 5, 7, 8, 11, {\bf 12}\}=2$\\  \hline 
& $e=[1 1 0 0 0 0]$ 	& No 		& Bad row \\  \hline
& $e'=[1 1 0 0 0 1]$ & No 		& $\lambda\{1, 2, 5, 7, 8, 11, {\bf 12}\}=2$\\  \hline
& $3'=[1 1 0 1 0 0]$  & Yes 	& - \\  \hline
& $4'=[1 1 1 1 0 0]$  & Yes 	& - \\  \hline
& $9'=[0 0 1 0 0 1]$  & Yes 	& - \\  \hline
& $10'=[0 0 0 1 0 1]$  & Yes 	& - \\  \hline \hline

$A$ with column $[0 0 1 1 0]$  & $b=[1 0 0 0 1 0]$ 	& No 		&  $\lambda\{1, 2, 5, {\bf 6}, 7, 8, 11 \}=2$	\\  \hline
and & $b'=[1 0 0 0 1 1]$ 	& Yes 		&  -	  \\ \hline
$B$ with column $[1 1 1 0 0]$ & $c=[1 1 0 0 1 0]$ 	& No 			&  	$\lambda\{1, 2, 5, {\bf 6}, 7, 8, 11 \}=2$	\\  \hline
& $c'=[1 1 0 0 1 1]$ & yes  		& -   		\\  \hline \hline

$C$ with column $[1 1 0 0 0]$  & $b=[1 0 0 0 1 0]$ 	& No 		&  $\lambda\{1, 2, 5, {\bf 6}, 7, 8, 11 \}=2$	\\  \hline
& $b'=[1 0 0 0 1 1]$ 	& No		&  Bad row	  \\ \hline
& $c=[1 1 0 0 1 0]$ 	& No 			&  	$\lambda\{1, 2, 5, {\bf 6}, 7, 8, 11 \}=2$	\\  \hline
& $c'=[1 1 0 0 1 1]$ & No 		& Bad row   		\\  \hline \hline
\end{tabular}
\end{center}
\normalsize
\small
 \begin{center} Table 2a: Cosimple single-element coextensions of $A_{22}$, $A_{11}$, $A_5$   in $EX[S_{10}$, $S^*_{10}]$ with respect to $(A_1, B_1)$ \end{center} 
\normalsize
\smallskip

\small
 \begin{center}
\begin{tabular}{|c|p{10em}|p{6em}|c|p{13em}|}
\hline
\bf   &\bf{Coextension Rows} &\bf{$S_{10}$, $S^*_{10}$ } & \bf {$A_1=\{1, 2, 3, 4, 9, 10\}$} \\  \hline \hline

$A$ with column $[0 0 1 1 0]$  & $a=[0 0 1 1 0 0]$  & No &   Bad row \\ \hline
and & $a'=[0 0 1 1 0 1]$ 	& No 	& Bad row         \\  \hline
$B$ with column $[1 1 1 0 0]$  & $b=[1 0 0 0 1 0]$ 	& No 		&  $\{1, 2, 3, 4, 9, 10, {\bf 12}\}$		\\  \hline
& $b'=[1 0 0 0 1 1]$ 	& Yes 		&  -		 \\  \hline
& $c=[1 1 0 0 1 0]$ 	& No 		& $\{1, 2, 3, 4, 9, 10, {\bf 12}\}$	 		\\  \hline
& $c'=[1 1 0 0 1 1]$ & yes 	& -   		\\  \hline
& $d=[1 1 1 0 0 0]$  & No 	& Bad row          \\  \hline  
& $d'=[1 1 1 0 0 1]$ & No 		& Bad row\\  \hline 
& $e=[1 1 0 0 0 0]$ 	& No 		& $\{1, 2, 3, 4, 9, 10, {\bf 12}\}$	 \\  \hline
& $e'=[1 1 0 0 0 1]$ & No 		& Bad row\\  \hline
& $3'=[1 1 0 1 0 0]$  & Yes 	& - \\  \hline
& $4'=[1 1 1 1 0 0]$  & Yes 	& - \\  \hline
& $9'=[0 0 1 0 0 1]$  & Yes 	& - \\  \hline
& $10'=[0 0 0 1 0 1]$  & Yes 	& - \\  \hline \hline

$A$ with column $[1 0 1 1 0]$  & $a=[0 0 1 1 0 0]$ 	& No 		&  $\{1, 2, 3, 4, {\bf 6}, 9, 10\}$	\\  \hline
and & $b'=[1 0 0 0 1 1]$ 	& Yes 		&  -	  \\ \hline
$B$ with column $[0 1 1 1 1]$& $c=[1 1 0 0 1 0]$ 	& Yes		   &  	-	\\  \hline
& $c'=[1 1 0 0 1 1]$ & No  		& $\{1, 2, 3, 4, {\bf 6}, 9, 10\}$   		\\  \hline \hline

$C$ with column $[1 1 0 0 0]$  & $a=[0 0 1 1 0 0]$ 	& No 		&  $\lambda\{1, 2, 5, {\bf 6}, 7, 8, 11 \}=2$	\\  \hline
& $a'=[0 0 1 1 0 1]$ 	& No		&  Bad row	  \\ \hline
& $d=[1 1 1 0 0 0]$ 	& No 			&  	$\lambda\{1, 2, 5, {\bf 6}, 7, 8, 11 \}=2$	\\  \hline
& $d'=[1 1 1 0 0 1]$ & No 		& Bad row   		\\  \hline \hline
\end{tabular}
\end{center}
\normalsize
\small
 \begin{center} Table 2b: Cosimple single-element coextensions of $A_{22}$, $A_{11}$, $A_5$   in $EX[S_{10}$, $S^*_{10}]$ with respect to $(A_2, B_2)$ \end{center} 
\normalsize

\smallskip
Let us call the rows that give matrices in $EX[S_{10}, S_{10}^*]$ satisfying  conditions in Theorem 2.1(iii) for a particular 3-separation $(A_i, B_i)$ {\it good rows for $(A_i, B_i)$} and those that don't, but are not bridging matroids, {\it bad rows for $(A_i, B_i)$}. Observe from Table 2a and 2b for 3-separation $(A_1, B_1)$, $A$ with column $[0 0 1 1 0 0]$ and $B$ with column $[0 1 1 1 1]$ have the property that all rows to be added are good. Similarly for 3-separation $(A_2, B_2)$, $A$ with column $[1 0 1 1 0]$ and $B$ with column $[0 1 1 1 1]$ have the property that all rows to be added are good. Moreover, for $C$ with column $[11000]$ the sets of bad rows for $(A_i, B_i)$ are disjoint. It follows from Corollary 2.2 that either $(A_1, B_1)$ or $(A_2, B_2)$ is induced in all matroids in $EX[S_{10}$, $S^*_{10}]$ with an $E_4$-minor. 

Lastly, although we do not need to check Theorem 2.1(iv) since $E_4$ is self-dual, we did verify it by repeating all the calculations for the simple single-element extensions of $A^*$, $B^*$, and $C^*$. This completes the proof of Theorem 1.1. $\qed$


\bigskip

 \bigskip

\noindent {\bf References}

\begin{enumerate}

 \item \label{Kingan1997} Kingan, S. R. (1997) A generalization of a graph result by D. W. Hall, {\it Discrete Mathematics} {\bf 173},  129-135. 

\item \label{Kingan2002}  Kingan, S. R. and Lemos, M. (2002) Almost-graphic matroids, {\it Advances in Applied Mathematics}, {\bf 28}, 438 - 477.

\item 	\label{Kingansubmitted} Kingan, S. R. (submitted). On Seymour's Decomposition Theorem.

\item 	\label{Mayhewsubmitted} Mayhew, D. Royle, G, and Whittle, G. (submitted) Excluding Kuratowski graphs and their duals from binary matroids.

\item \label{Oxley1987} Oxley, J. G. (1987) The binary matroids with no 4-wheel minor, {\it Trans. Amer. Math. Soc.}  {\bf 154}, 63-75.

\item \label{Oxley2012} Oxley, J. G. (2012) {\it Matroid Theory}, Second Edition, Oxford University Press, New York.

\item \label{Qin2004} Qin, H. and Zhou, X. (2004) The class of binary matroids with no $M(K_{3,3})$-,  $M^*(K_{3,3})$-, $M(K_5)$- or $M^*(K_5)$-minor. {\it  J. Comb. Theory, Ser. B} {\bf 90(1)}, 173-184

\item \label{Zhou2004} Zhou, X. (2004) On internally 4-connected non-regular binary matroids, {\it  J. Combin. Theory   Ser. B}  {\bf 91 }, 327-343.   

\end{enumerate}

\end {document}